\documentclass[11pt]{article} 

\usepackage{amsthm}

\usepackage{amsmath,amssymb,amsfonts,dsfont} 
\usepackage[usenames,dvipsnames]{color}
\usepackage{enumerate}
\usepackage{hyperref}

\newcommand{\R}{\mathbb{R}}

\newcommand{\N}{\mathbb{N}}

\newcommand{\xx}{\mathbf{x}}

\newcommand{\ud}{\,\mathrm{d}}

\newcommand{\Reff}{R_{\rm eff}}
\newcommand{\Rave}{R_{\rm ave}}

\newcommand{\Raveb}[1]{\Rave \l(#1\r)}

\newcommand{\Ravehydro}[1]{ \gamma\l(#1\r) }

\newcommand{\TMd}{T_{M^d}}

\newcommand{\TM}{T_M}

\newcommand{\TMdue}{ T_{M^2} }

\newcommand{\Rmax}{R_{\rm max}}

\newcommand{\ppii}{\mathbf{\pi}}
\newcommand{\regdeg}{\delta}

\newcommand{\HIT}{H}
\newcommand{\COMM}{C}

    \newcommand{\map}[3]{#1: #2 \rightarrow #3}

    \renewcommand{\l}{\left}
    \renewcommand{\r}{\right}

\title{Note on ``Average resistance of toroidal graphs'' by Rossi, Frasca and Fagnani%
\thanks{%
The authors are indebted to Alex Olshevsky for bringing Theorem 6.1 in~\cite{AKC-PR89} to their attention.}
}

\author{Wilbert Samuel Rossi%
\thanks{W. S. Rossi is with Department of Applied Mathematics, University of Twente, 7500 AE Enschede, The Netherlands;
{\tt\small w.s.rossi@utwente.nl}.}
\and Paolo Frasca%
\thanks{P. Frasca is with GIPSA-lab, Univ. Grenoble Alpes, CNRS, Inria, F-38000 Grenoble, France; {\tt\small paolo.frasca@gipsa-lab.fr}}%
\and Fabio Fagnani%
\thanks{F. Fagnani is with Dipartimento di Scienze Matematiche, 
Politecnico di Torino, corso Duca degli Abruzzi 24, 10129 Torino, Italy; 
{\tt\small fabio.fagnani@polito.it}.}%
}

\begin{document}
\maketitle

%
%
%
%

%

\begin{abstract}
\noindent In our recent paper {\it W.S. Rossi, P. Frasca and F. Fagnani, ``Average resistance of toroidal graphs'', {\em SIAM Journal on Control and Optimization}, 53(4):2541--2557, 2015}, we studied how the average resistances of $d$-dimensional toroidal grids depend on the graph topology and on the dimension of the graph. Our results were based on the connection between resistance and Laplacian eigenvalues.
In this note, we contextualize our work in the body of literature about random walks on graphs. Indeed, the average effective resistance of the $d$-dimensional toroidal grid is proportional to the mean hitting time of the simple random walk on that grid. If $d\geq3 $, then the average resistance can be bounded uniformly in the number of nodes and its value is of order $1/d$ for large $d$. 
%
\end{abstract}

%

%
%

\section{Introduction}
In this note we clarify the relation of~\cite{RFF:2015:average} with some literature in the field of probability and we show that Conjecture~2.2 stated therein follows directly from Theorem 6.1 in~\cite{AKC-PR89}, a result proved with electrical techniques.
We begin by briefly recalling some definitions and facts on electrical networks and on random walks. Next, we focus our discussion on toroidal grids.
\section{Electrical Networks and Simple Random Walks}
We consider an undirected connected graph $G=(V,E)$ with $N = |V|$ vertices and $|E|$ edges.
We denote with $d_v$ the degree of $v$ and we call $\regdeg$-regular a graph with $d_v = \delta$ for every $v\in V$.
We may think of the graph as an electrical network with all edges having unit resistance.
Given two distinct vertices $u,v\in V$, we denote by $\Reff(u,v)$ the effective resistance between them, i.e. the electrical potential difference induced between $u$ and $v$ by a unit current injected in $u$ and extracted from $v$. 
We define the average effective resistance of $G$ as 
\begin{align} \label{eq:rave-def}
	\Raveb{G}:=\frac{1}{2N^2}\sum\limits_{u,v\in V}\Reff(u,v)
\end{align}
where $\Reff(u,u):=0$. 
A \textit{simple random walk} on the graph $G = (V,E)$ is a  
discrete-time random process started at one vertex in $V$. 
See~\cite{LL:1993} for a concise introduction. 
If at step $t$ the random walk is at a node $v$, at step $t+1$ the walk will be in one of the neighbors of $v$ chosen with 
probability $1 / d_v$.
The sequence of vertices visited by the simple random walk is a time-reversible Markov chain with transition matrix 
\begin{align*}
	P_{v,w} = \l\{ \begin{array}{ll} d_v^{-1}  & \textup{if~} \{v,w\} \in E \\ 
			0 & \textup{if not}. \end{array} \r.
\end{align*}
Associated with a time-reversible Markov chain there is a unique positive, sum-one vector $\ppii$ called \textit{stationary distribution} satisfying the equations $\pi_v P_{v,w} = \pi_w P_{w,v}$ for all pairs in $V\times V$. 
For the simple random walk it is easy to show that $\pi_v = d_v({2|E|})^{-1} $. 
Given a pair of vertices $v,w \in V$, the  \textit{hitting time} $\HIT_{vw}$ is  the expected number of steps it takes to a random walk started in $v$ to first reach $w$. 
The sum $\COMM_{vw} = \HIT_{vw} + \HIT_{wv} $ is called \textit{commute time} and represents the expected number of steps that the random walk started at $v$ takes to reach $w$ and get back to $v$. 
For any pair of vertices $v,w$, the commute time $\COMM_{vw}$ and the effective resistance $\Reff(v,w)$ are proportional \cite[Thm.~2.1]{AKC-PR89}: 
$\COMM_{vw} = 2 |E| \Reff(v,w). $
%
Following  \cite[p. 113]{AF:94} we also introduce the \textit{average hitting time}
$$\tau_0(G) := \sum_{v,w} \pi_v\pi_w  \HIT_{vw}. $$
If the graph $G = (V,E)$ is $\regdeg$-regular, then the average hitting time and the average effective resistance are proportional.
Moreover, the stationary distribution $\ppii$ becomes uniform: 
$ \pi_v = N^{-1}, \forall v$.    
Hence, from the definition of $\tau_0$ 
\begin{align*}
	\tau_0(G)	&= \sum_{v,w} \pi_v\pi_w  \HIT_{vw}
			 = \frac{1}{2N^2} \sum_{v,w} \HIT_{vw} + \HIT_{wv}
			 = \frac{1}{2N^2} \sum_{v,w} \COMM_{vw}  \\
			& = \frac{1}{2N^2} \sum_{v,w}  2 |E| \Reff(v,w) 
		= \regdeg N \Raveb{G}. 
\end{align*}

\section{$d$-dimensional Toroidal Grids}  
%
%

The toroidal $d$-dimensional grids $\TMd$ is the graph where the vertices are arranged on a Cartesian lattice in $d$ dimensions with sides of length $M\geq3$ and periodic boundary conditions. Hence, $\TMd$ is $2d$-regular and contains $N = M^d$ vertices.
%
The equality 
\begin{align*} 
	\tau_0(\TMd) = 2 d M^d \Raveb{\TMd}
\end{align*}
allows us to interpret the results about $\tau_0(\TMd)$ in \cite[Prop. 13.8]{AF:94} as asymptotic relations\footnote{Given two sequences $\map{f,g}{\N}{\R^+}$, let $\ell^+=\limsup_n f(n)/g(n)$ and $\ell^-=\liminf_n f(n)/g(n)$. We write that $f=O(g)$ when $\ell^+<+\infty$; that  $f=o(g)$ when $\ell^+=0$; that $f \sim g$ when $\ell^+=\ell^-=1$, and $f=\Theta(g)$ when $\ell^+,\, \ell^-\in (0,+\infty)$.} for the average effective resistance of toroidal grids with $M$ growing large. 
The asymptotic trends are
\begin{alignat}{3}
	& \Raveb{\TM} 	&&	\sim \frac{1}{12} M 	&& \label{eq:AF-TM}	\\	
	& \Raveb{\TMdue}&&	\sim \frac{1}{ 2 \pi}\log M && \label{eq:AF-TM2}	\\  
	& \Raveb{\TMd} 	&&	\sim \Ravehydro{d}		&& \textup{~~~for~} d\geq 3  \label{eq:AF-TMd}
\end{alignat}
where the integral
\begin{align*} 
     \Ravehydro{d} := \int_{ \l[0, 1 \r]^d} \frac{1 }{2 d - 2\sum_{i=1}^d \cos(2\pi x_i)} \ud \xx
\end{align*}
converges for $d \geq 3$ with $2 d \Ravehydro{d} < \infty$.
%
%
%
%
These asymptotic trends are also obtained in our paper, by different techniques. The most interesting is~\eqref{eq:AF-TMd}, which in~\cite{AF:94} is directly proved as an asymptotic relation by a probabilistic argument. Instead, in~\cite{RFF:2015:average} we first prove that $\Ravehydro{d} \in [\frac{1}{4d}, \frac{4}{d} ]$ for $d \geq 3$ (Lem~3.1), and then we use these bounds to obtain estimates of $\Raveb{\TMd}$ for finite $M$ (Thm~2.4). 
Those estimates are such that 
$$\frac{1}{4d} \leq \lim_{M\to +\infty} \Raveb{\TMd} \leq \frac{8}{d+1} \qquad {\rm for~} d\geq 3 $$
and this relation asymptotically implies~\eqref{eq:AF-TMd}.
Additionally, we conjectured (Conjecture~2.2 of~\cite{RFF:2015:average}) that 
$$\Raveb{\TMd}=\Theta\left(\frac{1}{d}\right)\quad{\rm for}\; d\to +\infty\,,\; M\;{\rm fixed}.$$
Also this stronger statement can actually be proved, as follows.
For the lower bound, one can easily see~\cite[Thm. 2.4]{RFF:2015:average} that $ \Raveb{\TMd} \geq  \frac{1}{4d}$.
For the upper bound, let $\Rmax(G):= \max_{u,v\in V}\Reff(u,v)$ be the maximum of the effective resistances between any pair of vertices in the graph $G$ and observe that $\Reff(u,v) \leq \frac{1}{2} \Rmax(G)$.
Since Theorem~6.1 of~\cite{AKC-PR89} states that
$$\Rmax(\TMd) = \Theta\left({d}^{-1}\right) \quad{\rm d\geq3},$$ the conjecture is proved.

%
%

%
%
%
%
%
%
%
%

\bibliographystyle{plain}

\end{document}